\documentclass{amsart}
\usepackage{amsfonts}

\newtheorem{theorem}{Theorem}
\theoremstyle{plain}

\newtheorem{corollary}{Corollary}

\newtheorem{lemma}{Lemma}

\newtheorem{remark}{Remark}

\numberwithin{equation}{section}
\input{tcilatex}

\begin{document}
\title[Boas-Bellman Inequality]{On the Boas-Bellman Inequality in Inner Product Spaces}
\author{S.S. Dragomir}
\address{School of Computer Science and Mathematics\\
Victoria University of Technology\\
PO Box 14428, MCMC \\
Victoria 8001, Australia.}
\email{sever.dragomir@vu.edu.au}
\urladdr{http://rgmia.vu.edu.au/SSDragomirWeb.html}
\date{3 June, 2003.}

\begin{abstract}
New results related to the Boas-Bellman generalisation of Bessel's
inequality in inner product spaces are given.
\end{abstract}

\maketitle

\section{Introduction}

Let $\left( H;\left( \cdot ,\cdot \right) \right) $ be an inner product
space over the real or complex number field $\mathbb{K}$. If $\left(
e_{i}\right) _{1\leq i\leq n}$ are orthonormal vectors in the inner product
space $H,$ i.e., $\left( e_{i},e_{j}\right) =\delta _{ij}$ for all $i,j\in
\left\{ 1,\dots ,n\right\} $ where $\delta _{ij}$ is the Kronecker delta,
then the following inequality is well known in the literature as Bessel's
inequality (see for example \cite[p. 391]{6b}): 
\begin{equation}
\sum_{i=1}^{n}\left| \left( x,e_{i}\right) \right| ^{2}\leq \left\|
x\right\| ^{2}\text{ \ for any \ }x\in H.  \label{1.1}
\end{equation}

For other results related to Bessel's inequality, see \cite{3b} -- \cite{5b}
and Chapter XV in the book \cite{6b}.

In 1941, R.P. Boas \cite{2b} and in 1944, independently, R. Bellman \cite{1b}
proved the following generalisation of Bessel's inequality (see also 
\cite[p. 392]{6b}).

\begin{theorem}
\label{t1.1}If $x,y_{1},\dots ,y_{n}$ are elements of an inner product space 
$\left( H;\left( \cdot ,\cdot \right) \right) ,$ then the following
inequality: 
\begin{equation}
\sum_{i=1}^{n}\left| \left( x,y_{i}\right) \right| ^{2}\leq \left\|
x\right\| ^{2}\left[ \max_{1\leq i\leq n}\left\| y_{i}\right\| ^{2}+\left(
\sum_{1\leq i\neq j\leq n}\left| \left( y_{i},y_{j}\right) \right|
^{2}\right) ^{\frac{1}{2}}\right]  \label{1.2}
\end{equation}
holds.
\end{theorem}

A recent generalisation of the Boas-Bellman result was given in
Mitrinovi\'{c}-Pe\v{c}ari\'{c}-Fink \cite[p. 392]{6b} where they proved the
following.

\begin{theorem}
\label{t1.2}If $x,y_{1},\dots ,y_{n}$ are as in Theorem \ref{t1.1} and $%
c_{1},\dots ,c_{n}\in \mathbb{K}$, then one has the inequality: 
\begin{equation}
\left\vert \sum_{i=1}^{n}c_{i}\left( x,y_{i}\right) \right\vert ^{2}\leq
\left\Vert x\right\Vert ^{2}\sum_{i=1}^{n}\left\vert c_{i}\right\vert ^{2} 
\left[ \max_{1\leq i\leq n}\left\Vert y_{i}\right\Vert ^{2}+\left(
\sum_{1\leq i\neq j\leq n}\left\vert \left( y_{i},y_{j}\right) \right\vert
^{2}\right) ^{\frac{1}{2}}\right] .  \label{1.3}
\end{equation}
\end{theorem}

They also noted that if in (\ref{1.3}) one chooses $c_{i}=\overline{\left(
x,y_{i}\right) }$, then this inequality becomes (\ref{1.2}).

For other results related to the Boas-Bellman inequality, see \cite{4b}.

In this paper we point out some new results that may be related to both the
Mitrinovi\'{c}-Pe\v{c}ari\'{c}-Fink and Boas-Bellman inequalities.

\section{Some Preliminary Results}

We start with the following lemma which is also interesting in itself.

\begin{lemma}
\label{l2.1}Let $z_{1},\dots ,z_{n}\in H$ and $\alpha _{1},\dots ,\alpha
_{n}\in \mathbb{K}$. Then one has the inequality: 
\begin{multline}
\left\| \sum_{i=1}^{n}\alpha _{i}z_{i}\right\| ^{2}  \label{2.1} \\
\leq \left\{ 
\begin{array}{l}
\max\limits_{1\leq i\leq n}\left| \alpha _{i}\right|
^{2}\sum\limits_{i=1}^{n}\left\| z_{i}\right\| ^{2}; \\ 
\\ 
\left( \sum\limits_{i=1}^{n}\left| \alpha _{i}\right| ^{2\alpha }\right) ^{%
\frac{1}{\alpha }}\left( \sum\limits_{i=1}^{n}\left\| z_{i}\right\| ^{2\beta
}\right) ^{\frac{1}{\beta }},\ \ \ \text{where \ }\alpha >1,\frac{1}{\alpha }%
+\frac{1}{\beta }=1; \\ 
\\ 
\sum\limits_{i=1}^{n}\left| \alpha _{i}\right| ^{2}\max\limits_{1\leq i\leq
n}\left\| z_{i}\right\| ^{2},
\end{array}
\right. \\
+\left\{ 
\begin{array}{l}
\max\limits_{1\leq i\neq j\leq n}\left\{ \left| \alpha _{i}\alpha
_{j}\right| \right\} \sum\limits_{1\leq i\neq j\leq n}\left| \left(
z_{i},z_{j}\right) \right| ; \\ 
\\ 
\left[ \left( \sum\limits_{i=1}^{n}\left| \alpha _{i}\right| ^{\gamma
}\right) ^{2}-\left( \sum\limits_{i=1}^{n}\left| \alpha _{i}\right|
^{2\gamma }\right) \right] ^{\frac{1}{\gamma }}\left( \sum\limits_{1\leq
i\neq j\leq n}\left| \left( z_{i},z_{j}\right) \right| ^{\delta }\right) ^{%
\frac{1}{\delta }}, \\ 
\hfill \ \ \ \text{where \ }\gamma >1,\ \ \frac{1}{\gamma }+\frac{1}{\delta }%
=1; \\ 
\\ 
\left[ \left( \sum\limits_{i=1}^{n}\left| \alpha _{i}\right| \right)
^{2}-\sum\limits_{i=1}^{n}\left| \alpha _{i}\right| ^{2}\right]
\max\limits_{1\leq i\neq j\leq n}\left| \left( z_{i},z_{j}\right) \right| .
\end{array}
\right.
\end{multline}
\end{lemma}

\begin{proof}
We observe that 
\begin{align}
\left\| \sum_{i=1}^{n}\alpha _{i}z_{i}\right\| ^{2}& =\left(
\sum_{i=1}^{n}\alpha _{i}z_{i},\sum_{j=1}^{n}\alpha _{j}z_{j}\right)
\label{2.2} \\
& =\sum_{i=1}^{n}\sum_{j=1}^{n}\alpha _{i}\overline{\alpha _{j}}\left(
z_{i},z_{j}\right) =\left| \sum_{i=1}^{n}\sum_{j=1}^{n}\alpha _{i}\overline{%
\alpha _{j}}\left( z_{i},z_{j}\right) \right|  \notag \\
& \leq \sum_{i=1}^{n}\sum_{j=1}^{n}\left| \alpha _{i}\right| \left| 
\overline{\alpha _{j}}\right| \left| \left( z_{i},z_{j}\right) \right| 
\notag \\
& =\sum_{i=1}^{n}\left| \alpha _{i}\right| ^{2}\left\| z_{i}\right\|
^{2}+\sum\limits_{1\leq i\neq j\leq n}\left| \alpha _{i}\right| \left|
\alpha _{j}\right| \left| \left( z_{i},z_{j}\right) \right| .  \notag
\end{align}
Using H\"{o}lder's inequality, we may write that 
\begin{eqnarray}
&&\sum_{i=1}^{n}\left| \alpha _{i}\right| ^{2}\left\| z_{i}\right\| ^{2}
\label{2.3} \\
&\leq &\left\{ 
\begin{array}{l}
\max\limits_{1\leq i\leq n}\left| \alpha _{i}\right|
^{2}\sum\limits_{i=1}^{n}\left\| z_{i}\right\| ^{2}; \\ 
\\ 
\left( \sum\limits_{i=1}^{n}\left| \alpha _{i}\right| ^{2\alpha }\right) ^{%
\frac{1}{\alpha }}\left( \sum\limits_{i=1}^{n}\left\| z_{i}\right\| ^{2\beta
}\right) ^{\frac{1}{\beta }},\ \ \ \text{where \ }\alpha >1,\frac{1}{\alpha }%
+\frac{1}{\beta }=1; \\ 
\\ 
\sum\limits_{i=1}^{n}\left| \alpha _{i}\right| ^{2}\max\limits_{1\leq i\leq
n}\left\| z_{i}\right\| ^{2}.
\end{array}
\right.  \notag
\end{eqnarray}
By H\"{o}lder's inequality for double sums we also have 
\begin{equation}
\sum\limits_{1\leq i\neq j\leq n}\left| \alpha _{i}\right| \left| \alpha
_{j}\right| \left| \left( z_{i},z_{j}\right) \right|  \label{2.4}
\end{equation}
\begin{eqnarray*}
&\leq &\left\{ 
\begin{array}{l}
\max\limits_{1\leq i\neq j\leq n}\left| \alpha _{i}\alpha _{j}\right|
\sum\limits_{1\leq i\neq j\leq n}\left| \left( z_{i},z_{j}\right) \right| ;
\\ 
\\ 
\left( \sum\limits_{1\leq i\neq j\leq n}\left| \alpha _{i}\right| ^{\gamma
}\left| \alpha _{j}\right| ^{\gamma }\right) ^{\frac{1}{\gamma }}\left(
\sum\limits_{1\leq i\neq j\leq n}\left| \left( z_{i},z_{j}\right) \right|
^{\delta }\right) ^{\frac{1}{\delta }}, \\ 
\hfill \ \ \ \text{where \ }\gamma >1,\ \ \frac{1}{\gamma }+\frac{1}{\delta }%
=1; \\ 
\\ 
\sum\limits_{1\leq i\neq j\leq n}\left| \alpha _{i}\right| \left| \alpha
_{j}\right| \max\limits_{1\leq i\neq j\leq n}\left| \left(
z_{i},z_{j}\right) \right| ,
\end{array}
\right. \\
&& \\
&=&\left\{ 
\begin{array}{l}
\max\limits_{1\leq i\neq j\leq n}\left\{ \left| \alpha _{i}\alpha
_{j}\right| \right\} \sum\limits_{1\leq i\neq j\leq n}\left| \left(
z_{i},z_{j}\right) \right| ; \\ 
\\ 
\left[ \left( \sum\limits_{i=1}^{n}\left| \alpha _{i}\right| ^{\gamma
}\right) ^{2}-\left( \sum\limits_{i=1}^{n}\left| \alpha _{i}\right|
^{2\gamma }\right) \right] ^{\frac{1}{\gamma }}\left( \sum\limits_{1\leq
i\neq j\leq n}\left| \left( z_{i},z_{j}\right) \right| ^{\delta }\right) ^{%
\frac{1}{\delta }}, \\ 
\hfill \ \ \ \text{where \ }\gamma >1,\ \ \frac{1}{\gamma }+\frac{1}{\delta }%
=1; \\ 
\\ 
\left[ \left( \sum\limits_{i=1}^{n}\left| \alpha _{i}\right| \right)
^{2}-\sum\limits_{i=1}^{n}\left| \alpha _{i}\right| ^{2}\right]
\max\limits_{1\leq i\neq j\leq n}\left| \left( z_{i},z_{j}\right) \right| .
\end{array}
\right.
\end{eqnarray*}
Utilising (\ref{2.3}) and (\ref{2.4}) in (\ref{2.2}), we may deduce the
desired result (\ref{2.1}).
\end{proof}

\begin{remark}
\label{r2.2}Inequality (\ref{2.1}) contains in fact 9 different inequalities
which may be obtained combining the first 3 ones with the last 3 ones.
\end{remark}

A particular case that may be related to the Boas-Bellman result is embodied
in the following inequality.

\begin{corollary}
\label{c2.3}With the assumptions in Lemma \ref{l2.1}, we have 
\begin{equation}
\left\Vert \sum_{i=1}^{n}\alpha _{i}z_{i}\right\Vert ^{2}  \label{2.5}
\end{equation}
\end{corollary}

\begin{align}
& \leq \sum_{i=1}^{n}\left| \alpha _{i}\right| ^{2}\left\{
\max\limits_{1\leq i\leq n}\left\| z_{i}\right\| ^{2}+\frac{\left[ \left(
\sum_{i=1}^{n}\left| \alpha _{i}\right| ^{2}\right)
^{2}-\sum_{i=1}^{n}\left| \alpha _{i}\right| ^{4}\right] ^{\frac{1}{2}}}{%
\sum_{i=1}^{n}\left| \alpha _{i}\right| ^{2}}\left( \sum\limits_{1\leq i\neq
j\leq n}\left| \left( z_{i},z_{j}\right) \right| ^{2}\right) ^{\frac{1}{2}%
}\right\}  \notag \\
& \leq \sum_{i=1}^{n}\left| \alpha _{i}\right| ^{2}\left\{
\max\limits_{1\leq i\leq n}\left\| z_{i}\right\| ^{2}+\left(
\sum\limits_{1\leq i\neq j\leq n}\left| \left( z_{i},z_{j}\right) \right|
^{2}\right) ^{\frac{1}{2}}\right\} .  \notag
\end{align}
The first inequality follows by taking the third branch in the first curly
bracket with the second branch in the second curly bracket for $\gamma
=\delta =2.$

The second inequality in (\ref{2.5}) follows by the fact that 
\begin{equation*}
\left[ \left( \sum_{i=1}^{n}\left| \alpha _{i}\right| ^{2}\right)
^{2}-\sum_{i=1}^{n}\left| \alpha _{i}\right| ^{4}\right] ^{\frac{1}{2}}\leq
\sum_{i=1}^{n}\left| \alpha _{i}\right| ^{2}.
\end{equation*}
Applying the following Cauchy-Bunyakovsky-Schwarz type inequality 
\begin{equation}
\left( \sum_{i=1}^{n}a_{i}\right) ^{2}\leq n\sum_{i=1}^{n}a_{i}^{2},\ \ \
a_{i}\in \mathbb{R}_{+},\ \ 1\leq i\leq n,  \label{2.6}
\end{equation}
we may write that 
\begin{equation}
\left( \sum\limits_{i=1}^{n}\left| \alpha _{i}\right| ^{\gamma }\right)
^{2}-\sum\limits_{i=1}^{n}\left| \alpha _{i}\right| ^{2\gamma }\leq \left(
n-1\right) \sum\limits_{i=1}^{n}\left| \alpha _{i}\right| ^{2\gamma }\ \ \ \
\ \ \left( n\geq 1\right)  \label{2.7}
\end{equation}
and 
\begin{equation}
\left( \sum\limits_{i=1}^{n}\left| \alpha _{i}\right| \right)
^{2}-\sum\limits_{i=1}^{n}\left| \alpha _{i}\right| ^{2}\leq \left(
n-1\right) \sum\limits_{i=1}^{n}\left| \alpha _{i}\right| ^{2}\ \ \ \ \ \
\left( n\geq 1\right) .  \label{2.8}
\end{equation}
Also, it is obvious that: 
\begin{equation}
\max\limits_{1\leq i\neq j\leq n}\left\{ \left| \alpha _{i}\alpha
_{j}\right| \right\} \leq \max\limits_{1\leq i\leq n}\left| \alpha
_{i}\right| ^{2}.  \label{2.9}
\end{equation}
Consequently, we may state the following coarser upper bounds for $\left\|
\sum_{i=1}^{n}\alpha _{i}z_{i}\right\| ^{2}$ that may be useful in
applications.

\begin{corollary}
\label{c2.4}With the assumptions in Lemma \ref{l2.1}, we have the
inequalities: 
\begin{multline}
\left\| \sum_{i=1}^{n}\alpha _{i}z_{i}\right\| ^{2}  \label{2.10} \\
\leq \left\{ 
\begin{array}{l}
\max\limits_{1\leq i\leq n}\left| \alpha _{i}\right|
^{2}\sum\limits_{i=1}^{n}\left\| z_{i}\right\| ^{2}; \\ 
\\ 
\left( \sum\limits_{i=1}^{n}\left| \alpha _{i}\right| ^{2\alpha }\right) ^{%
\frac{1}{\alpha }}\left( \sum\limits_{i=1}^{n}\left\| z_{i}\right\| ^{2\beta
}\right) ^{\frac{1}{\beta }},\ \ \ \text{where \ }\alpha >1,\frac{1}{\alpha }%
+\frac{1}{\beta }=1; \\ 
\\ 
\sum\limits_{i=1}^{n}\left| \alpha _{i}\right| ^{2}\max\limits_{1\leq i\leq
n}\left\| z_{i}\right\| ^{2},
\end{array}
\right. \\
+\left\{ 
\begin{array}{l}
\max\limits_{1\leq i\leq n}\left| \alpha _{i}\right| ^{2}\sum\limits_{1\leq
i\neq j\leq n}\left| \left( z_{i},z_{j}\right) \right| ; \\ 
\\ 
\left( n-1\right) ^{\frac{1}{\gamma }}\left( \sum\limits_{i=1}^{n}\left|
\alpha _{i}\right| ^{2\gamma }\right) ^{\frac{1}{\gamma }}\left(
\sum\limits_{1\leq i\neq j\leq n}\left| \left( z_{i},z_{j}\right) \right|
^{\delta }\right) ^{\frac{1}{\delta }}, \\ 
\hfill \ \ \ \text{where \ }\gamma >1,\ \ \frac{1}{\gamma }+\frac{1}{\delta }%
=1; \\ 
\\ 
\left( n-1\right) \sum\limits_{i=1}^{n}\left| \alpha _{i}\right|
^{2}\max\limits_{1\leq i\neq j\leq n}\left| \left( z_{i},z_{j}\right)
\right| .
\end{array}
\right.
\end{multline}
\end{corollary}

The proof is obvious by Lemma \ref{l2.1} in applying the inequalities (\ref
{2.7}) -- (\ref{2.9}).

\begin{remark}
\label{r2.5}The following inequalities which are incorporated in (\ref{2.10}%
) are of special interest: 
\begin{equation}
\left\| \sum_{i=1}^{n}\alpha _{i}z_{i}\right\| ^{2}\leq \max\limits_{1\leq
i\leq n}\left| \alpha _{i}\right| ^{2}\left[ \sum\limits_{i=1}^{n}\left\|
z_{i}\right\| ^{2}+\sum\limits_{1\leq i\neq j\leq n}\left| \left(
z_{i},z_{j}\right) \right| \right] ;  \label{2.11}
\end{equation}
\begin{multline}
\left\| \sum_{i=1}^{n}\alpha _{i}z_{i}\right\| ^{2}  \label{2.12} \\
\leq \left( \sum\limits_{i=1}^{n}\left| \alpha _{i}\right| ^{2p}\right) ^{%
\frac{1}{p}}\left[ \left( \sum\limits_{i=1}^{n}\left\| z_{i}\right\|
^{2q}\right) ^{\frac{1}{q}}+\left( n-1\right) ^{\frac{1}{p}}\left(
\sum\limits_{1\leq i\neq j\leq n}\left| \left( z_{i},z_{j}\right) \right|
^{q}\right) ^{\frac{1}{q}}\right] ,
\end{multline}
where $p>1,$ $\frac{1}{p}+\frac{1}{q}=1;$ and 
\begin{equation}
\left\| \sum_{i=1}^{n}\alpha _{i}z_{i}\right\| ^{2}\leq
\sum\limits_{i=1}^{n}\left| \alpha _{i}\right| ^{2}\left[ \max\limits_{1\leq
i\leq n}\left\| z_{i}\right\| ^{2}+\left( n-1\right) \max\limits_{1\leq
i\neq j\leq n}\left| \left( z_{i},z_{j}\right) \right| \right] .
\label{2.13}
\end{equation}
\end{remark}

\section{Some Mitrinovi\'{c}-Pe\v{c}ari\'{c}-Fink Type Inequalities}

We are now able to point out the following result which complements the
inequality (\ref{1.3}) due to Mitrinovi\'{c}, Pe\v{c}ari\'{c} and Fink 
\cite[p. 392]{6b}.

\begin{theorem}
\label{t3.1}Let $x,y_{1},\dots ,y_{n}$ be vectors of an inner product space $%
\left( H;\left( \cdot ,\cdot \right) \right) $ and $c_{1},\dots ,c_{n}\in 
\mathbb{K}$ $\left( \mathbb{K}=\mathbb{C},\mathbb{R}\right) .$ Then one has
the inequalities: 
\begin{multline}
\left| \sum_{i=1}^{n}c_{i}\left( x,y_{i}\right) \right| ^{2}  \label{3.1} \\
\leq \left\| x\right\| ^{2}\times \left\{ 
\begin{array}{l}
\max\limits_{1\leq i\leq n}\left| c_{i}\right|
^{2}\sum\limits_{i=1}^{n}\left\| y_{i}\right\| ^{2}; \\ 
\\ 
\left( \sum\limits_{i=1}^{n}\left| c_{i}\right| ^{2\alpha }\right) ^{\frac{1%
}{\alpha }}\left( \sum\limits_{i=1}^{n}\left\| y_{i}\right\| ^{2\beta
}\right) ^{\frac{1}{\beta }},\ \ \ \text{where \ }\alpha >1,\frac{1}{\alpha }%
+\frac{1}{\beta }=1; \\ 
\\ 
\sum\limits_{i=1}^{n}\left| c_{i}\right| ^{2}\max\limits_{1\leq i\leq
n}\left\| y_{i}\right\| ^{2},
\end{array}
\right. \\
+\left\| x\right\| ^{2}\times \left\{ 
\begin{array}{l}
\max\limits_{1\leq i\neq j\leq n}\left\{ \left| c_{i}c_{j}\right| \right\}
\sum\limits_{1\leq i\neq j\leq n}\left| \left( y_{i},y_{j}\right) \right| ;
\\ 
\\ 
\left[ \left( \sum\limits_{i=1}^{n}\left| c_{i}\right| ^{\gamma }\right)
^{2}-\left( \sum\limits_{i=1}^{n}\left| c_{i}\right| ^{2\gamma }\right) 
\right] ^{\frac{1}{\gamma }}\left( \sum\limits_{1\leq i\neq j\leq n}\left|
\left( y_{i},y_{j}\right) \right| ^{\delta }\right) ^{\frac{1}{\delta }}, \\ 
\hfill \ \ \ \text{where \ }\gamma >1,\ \ \frac{1}{\gamma }+\frac{1}{\delta }%
=1; \\ 
\\ 
\left[ \left( \sum\limits_{i=1}^{n}\left| c_{i}\right| \right)
^{2}-\sum\limits_{i=1}^{n}\left| c_{i}\right| ^{2}\right] \max\limits_{1\leq
i\neq j\leq n}\left| \left( y_{i},y_{j}\right) \right| .
\end{array}
\right.
\end{multline}
\end{theorem}

\begin{proof}
We note that 
\begin{equation*}
\sum_{i=1}^{n}c_{i}\left( x,y_{i}\right) =\left( x,\sum_{i=1}^{n}\overline{%
c_{i}}y_{i}\right) .
\end{equation*}
Using Schwarz's inequality in inner product spaces, we have: 
\begin{equation*}
\left| \sum_{i=1}^{n}c_{i}\left( x,y_{i}\right) \right| ^{2}\leq \left\|
x\right\| ^{2}\left\| \sum_{i=1}^{n}\overline{c_{i}}y_{i}\right\| ^{2}.
\end{equation*}
Now using Lemma \ref{l2.1} with $\alpha _{i}=\overline{c_{i}},$ $z_{i}=y_{i}$
$\left( i=1,\dots ,n\right) ,$ we deduce the desired inequality (\ref{3.1}).
\end{proof}

The following particular inequalities that may be obtained by the
Corollaries \ref{c2.3} and \ref{c2.4} and Remark \ref{r2.5} hold.

\begin{corollary}
\label{c3.2}With the assumptions in Theorem \ref{t3.1}, one has the
inequalities: 
\begin{multline}
\left| \sum_{i=1}^{n}c_{i}\left( x,y_{i}\right) \right| ^{2}  \label{3.2} \\
\leq \times \left\{ 
\begin{array}{l}
\left\| x\right\| ^{2}\sum\limits_{i=1}^{n}\left| c_{i}\right| ^{2}\left\{
\max\limits_{1\leq i\leq n}\left\| y_{i}\right\| ^{2}+\left(
\sum\limits_{1\leq i\neq j\leq n}\left| \left( y_{i},y_{j}\right) \right|
^{2}\right) ^{\frac{1}{2}}\right\} ; \\ 
\\ 
\left\| x\right\| ^{2}\max\limits_{1\leq i\leq n}\left| c_{i}\right|
^{2}\left\{ \sum\limits_{i=1}^{n}\left\| y_{i}\right\|
^{2}+\sum\limits_{1\leq i\neq j\leq n}\left| \left( y_{i},y_{j}\right)
\right| \right\} \\ 
\\ 
\left\| x\right\| ^{2}\left( \sum\limits_{i=1}^{n}\left| c_{i}\right|
^{2p}\right) ^{\frac{1}{p}}\left\{ \left( \sum\limits_{i=1}^{n}\left\|
y_{i}\right\| ^{2q}\right) ^{\frac{1}{q}}+\left( n-1\right) ^{\frac{1}{p}%
}\left( \sum\limits_{1\leq i\neq j\leq n}\left| \left( y_{i},y_{j}\right)
\right| ^{q}\right) ^{\frac{1}{q}}\right\} , \\ 
\hfill \ \ \ \text{where \ }p>1,\frac{1}{p}+\frac{1}{q}=1; \\ 
\left\| x\right\| ^{2}\sum\limits_{i=1}^{n}\left| c_{i}\right| ^{2}\left\{
\max\limits_{1\leq i\leq n}\left\| y_{i}\right\| ^{2}+\left( n-1\right)
\max\limits_{1\leq i\neq j\leq n}\left| \left( y_{i},y_{j}\right) \right|
\right\} ,
\end{array}
\right.
\end{multline}
\end{corollary}

\begin{remark}
\label{r3.3}Note that the first inequality in (\ref{3.2}) is the result
obtained by Mitrinovi\'{c}-Pe\v{c}ari\'{c}-Fink in \cite{6b}. The other 3
provide similar bounds in terms of the $p-$norms of the vector $\left(
\left| c_{1}\right| ^{2},\dots ,\left| c_{n}\right| ^{2}\right) .$
\end{remark}

\section{Some Boas-Bellman Type Inequalities}

If one chooses $c_{i}=\overline{\left( x,y_{i}\right) }$ $\left( i=1,\dots
,n\right) $ in (\ref{3.1}), then it is possible to obtain 9 different
inequalities between the Fourier coefficients $\left( x,y_{i}\right) $ and
the norms and inner products of the vectors $y_{i}$ $\left( i=1,\dots
,n\right) .$ We restrict ourselves only to those inequalities that may be
obtained from (\ref{3.2}).

As Mitrinovi\'{c}, Pe\v{c}ari\'{c} and Fink noted in \cite[p. 392]{6b}, the
first inequality in (\ref{3.2}) for the above selection of $c_{i}$ will
produce the Boas-Bellman inequality (\ref{1.2}).

From the second inequality in (\ref{3.2}) for $c_{i}=\overline{\left(
x,y_{i}\right) }$ we get 
\begin{equation*}
\left( \sum_{i=1}^{n}\left| \left( x,y_{i}\right) \right| ^{2}\right)
^{2}\leq \left\| x\right\| ^{2}\max_{1\leq i\leq n}\left| \left(
x,y_{i}\right) \right| ^{2}\left\{ \sum_{i=1}^{n}\left\| y_{i}\right\|
^{2}+\sum_{1\leq i\neq j\leq n}\left| \left( y_{i},y_{j}\right) \right|
\right\} .
\end{equation*}
Taking the square root in this inequality we obtain: 
\begin{equation}
\sum_{i=1}^{n}\left| \left( x,y_{i}\right) \right| ^{2}\leq \left\|
x\right\| \max_{1\leq i\leq n}\left| \left( x,y_{i}\right) \right| \left\{
\sum_{i=1}^{n}\left\| y_{i}\right\| ^{2}+\sum_{1\leq i\neq j\leq n}\left|
\left( y_{i},y_{j}\right) \right| \right\} ^{\frac{1}{2}},  \label{4.1}
\end{equation}
for any $x,y_{1},\dots ,y_{n}$ vectors in the inner product space $\left(
H;\left( \cdot ,\cdot \right) \right) .$

If we assume that $\left( e_{i}\right) _{1\leq i\leq n}$ is an orthonormal
family in $H,$ then by (\ref{4.1}) we have 
\begin{equation}
\sum_{i=1}^{n}\left| \left( x,e_{i}\right) \right| ^{2}\leq \sqrt{n}\left\|
x\right\| \max_{1\leq i\leq n}\left| \left( x,e_{i}\right) \right| ,\ \ \
x\in H.  \label{4.2}
\end{equation}
From the third inequality in (\ref{3.2}) for $c_{i}=\overline{\left(
x,y_{i}\right) }$ we deduce 
\begin{multline*}
\left( \sum_{i=1}^{n}\left| \left( x,y_{i}\right) \right| ^{2}\right)
^{2}\leq \left\| x\right\| ^{2}\left( \sum_{i=1}^{n}\left| \left(
x,y_{i}\right) \right| ^{2p}\right) ^{\frac{1}{p}} \\
\times \left\{ \left( \sum\limits_{i=1}^{n}\left\| y_{i}\right\|
^{2q}\right) ^{\frac{1}{q}}+\left( n-1\right) ^{\frac{1}{p}}\left(
\sum\limits_{1\leq i\neq j\leq n}\left| \left( y_{i},y_{j}\right) \right|
^{q}\right) ^{\frac{1}{q}}\right\} ,
\end{multline*}
for $p>1,$ $\frac{1}{p}+\frac{1}{q}=1.$

Taking the square root in this inequality we get 
\begin{multline}
\sum_{i=1}^{n}\left| \left( x,y_{i}\right) \right| ^{2}\leq \left\|
x\right\| \left( \sum_{i=1}^{n}\left| \left( x,y_{i}\right) \right|
^{2p}\right) ^{\frac{1}{2p}}  \label{4.3} \\
\times \left\{ \left( \sum\limits_{i=1}^{n}\left\| y_{i}\right\|
^{2q}\right) ^{\frac{1}{q}}+\left( n-1\right) ^{\frac{1}{p}}\left(
\sum\limits_{1\leq i\neq j\leq n}\left| \left( y_{i},y_{j}\right) \right|
^{q}\right) ^{\frac{1}{q}}\right\} ^{\frac{1}{2}},
\end{multline}
for any $x,y_{1},\dots ,y_{n}\in H,$ $p>1,$ $\frac{1}{p}+\frac{1}{q}=1.$

The above inequality (\ref{4.3}) becomes, for an orthornormal family $\left(
e_{i}\right) _{1\leq i\leq n},$%
\begin{equation}
\sum_{i=1}^{n}\left| \left( x,e_{i}\right) \right| ^{2}\leq n^{\frac{1}{q}%
}\left\| x\right\| \left( \sum_{i=1}^{n}\left| \left( x,e_{i}\right) \right|
^{2p}\right) ^{\frac{1}{2p}},\ \ \ x\in H.  \label{4.4}
\end{equation}
Finally, the choice $c_{i}=\overline{\left( x,y_{i}\right) }$ $\left(
i=1,\dots ,n\right) $ will produce in the last inequality in (\ref{3.2}) 
\begin{equation*}
\left( \sum_{i=1}^{n}\left| \left( x,y_{i}\right) \right| ^{2}\right)
^{2}\leq \left\| x\right\| ^{2}\sum_{i=1}^{n}\left| \left( x,y_{i}\right)
\right| ^{2}\left\{ \max\limits_{1\leq i\leq n}\left\| y_{i}\right\|
^{2}+\left( n-1\right) \max\limits_{1\leq i\neq j\leq n}\left| \left(
y_{i},y_{j}\right) \right| \right\}
\end{equation*}
giving the following Boas-Bellman type inequality 
\begin{equation}
\sum_{i=1}^{n}\left| \left( x,y_{i}\right) \right| ^{2}\leq \left\|
x\right\| ^{2}\left\{ \max\limits_{1\leq i\leq n}\left\| y_{i}\right\|
^{2}+\left( n-1\right) \max\limits_{1\leq i\neq j\leq n}\left| \left(
y_{i},y_{j}\right) \right| \right\} ,  \label{4.5}
\end{equation}
for any $x,y_{1},\dots ,y_{n}\in H.$

It is obvious that (\ref{4.5}) will give for orthonormal families the well
known Bessel inequality.

\begin{remark}
In order the compare the Boas-Bellman result with our result (\ref{4.5}), it
is enough to compare the quantities 
\begin{equation*}
A:=\left( \sum\limits_{1\leq i\neq j\leq n}\left| \left( y_{i},y_{j}\right)
\right| ^{2}\right) ^{\frac{1}{2}}
\end{equation*}
and 
\begin{equation*}
B:=\left( n-1\right) \max\limits_{1\leq i\neq j\leq n}\left| \left(
y_{i},y_{j}\right) \right| .
\end{equation*}
Consider the inner product space $H=\mathbb{R}$ with $\left( x,y\right) =xy,$
and choose $n=3,$ $y_{1}=a>0$, $y_{2}=b>0,$ $y_{3}=c>0.$ Then 
\begin{equation*}
A=\sqrt{2}\left( a^{2}b^{2}+b^{2}c^{2}+c^{2}a^{2}\right) ^{\frac{1}{2}},\ \
\ \ \ \ \ B=2\max \left( ab,ac,bc\right) .
\end{equation*}
Denote $ab=p,$ $bc=q,$ $ca=r.$ Then 
\begin{equation*}
A=\sqrt{2}\left( p^{2}+q^{2}+r^{2}\right) ^{\frac{1}{2}},\ \ \ \ \ \ \
B=2\max \left( p,q,r\right) .
\end{equation*}
Firstly, if we assume that $p=q=r,$ then $A=\sqrt{6}p,$ $B=2p$ which shows
that $A>B.$

Now choose $r=1$ and $p,q=\frac{1}{2}.$ Then $A=\sqrt{3}$ and $B=2$ showing
that $B>A.$

Consequently, in general, the Boas-Bellman inequality and our inequality (%
\ref{4.5}) cannot be compared.
\end{remark}

\end{document}